\documentclass[12pt]{article}
\usepackage{color}
\usepackage{amsmath}
\usepackage{latexsym}
\usepackage{amsfonts}
\usepackage{soul}
\usepackage{amssymb}

\newcommand{\scal}[1]{\langle#1\rangle}

\newcommand{\va}{\varepsilon}

\newtheorem{Pa}{Paper}[section]
\newtheorem{Tm}[Pa]{{\bf Theorem}}
\newtheorem{La}[Pa]{{\bf Lemma}}
\newtheorem{Cy}[Pa]{{\bf Corollary}}

\newtheorem{Pn}[Pa]{{\bf Proposition}}

\newtheorem{Pb}[Pa]{{\bf Problem}}

\newtheorem{Ex}[Pa]{{\bf Example}}

\newcommand{\mZ}{{\mathbb Z}}
\newcommand{\mC}{{\mathbb C}}

\newcommand{\mR}{{\mathbb R}}
\newcommand{\mT}{{\mathbb T}}

\newcommand{\qed}{\ifmmode$\Box$\else{\unskip\nobreak\hfil

\penalty50\hskip1em\null\nobreak\hfil$\Box$
\parfillskip=0pt\finalhyphendemerits=0\endgraf}\fi}

\newcommand{\gdots}{\mathinner{\mkern1mu\vbox{\kern7pt\hbox{.}}
\mkern2mu\raise2.5pt\hbox{.}\mkern2mu\raise5pt\hbox{.}\mkern1mu}}

\newcommand{\eq} [1] {\begin{equation}\label{#1}}
\newcommand{\en} {\end{equation}}
\newcommand {\eqn}      {\begin{eqnarray}}
\newcommand {\enn}      {\end{eqnarray}}
\newcommand {\bstar}    {\begin{eqnarray*}}
\newcommand {\estar}    {\end{eqnarray*}}
\newcommand {\mat}      [1] {\left[\begin{array}{#1}}
\newcommand {\rix}          {\end{array}\right]}

\begin{document}

\title{Non-denseness of factorable matrix functions}

\author{Alex Brudnyi\footnote{Research of this author is supported in part by NSERC.}
\\ Department of Mathematics \\ University of Calgary
\\ 2500 University Dr. NW \\ Calgary, Alberta, Canada T2N 1N4 \\ albru@math.ucalgary.ca
 \and Leiba Rodman \footnote{Research of this author supported in part by
the Faculty Research Assignment and Plumeri Award at the College of William and Mary.} \\
Department of Mathematics  \\
College of William and Mary\\
Williamsburg, VA 23187-8795,  USA \\ e-mail: lxrodm@math.wm.edu
\and Ilya M. Spitkovsky  \\Department of Mathematics  \\
College of William and Mary\\
Williamsburg, VA 23187-8795,  USA \\  e-mail:  ilya@math.wm.edu
}
\date{}

\maketitle

\begin{abstract}
It is proved that for certain algebras of continuous functions
on compact abelian groups, the set of factorable matrix functions with entries in the algebra is not dense in the group
of invertible matrix functions with entries in the algebra,
assuming  that the dual abelian group contains a subgroup isomorphic to ${\mathbb
Z}^3$. These algebras include the algebra of all continuous functions and
the Wiener algebra. More precisely, it is shown that infinitely many connected components
of the group of invertible matrix functions do not contain any factorable matrix functions,
again under the same assumption. Moreover, these components actually are disjoint with
the subgroup
generated by the triangularizable matrix functions.
\bigskip

\medskip

\end{abstract}

\section{Introduction}\label{jul51}
\setcounter{equation}{0}

Let $G$ be a (multiplicative) connected compact abelian group and let
$\Gamma$ be its (additive) character group. Recall that $\Gamma$
consists of continuous homomorphisms of $G$ into the group $\mT$ of
unimodular complex numbers. Note that $\Gamma$ is torsion free and discrete in the natural topology of
the dual locally compact group, and conversely every torsion free discrete abelian group is the character group of come
connected compact abelian group; see, e.g. \cite[Section 1.2]{Ru90}.

It is well-known \cite{Ru90} that, because $G$ is
connected, $\Gamma$ can be made into a linearly ordered group. So
let $\preceq$ be a fixed linear order such that $(\Gamma,\preceq)$ is
an ordered group.  Let $\Gamma_+=\{x\in \Gamma \, :\, x \succeq 0\}$,
$\Gamma_-=\{x\in \Gamma \, :\, x \preceq 0\}$. Standard widely used examples of $\Gamma$ are
${\mathbb Z}$ (the group of integers), ${\mathbb R}$ (the group of reals with the discrete topology),
and ${\mathbb Z}^k$, ${\mathbb R}^k$ with a lexicographic ordering (where $k$ is positive integer).

If $U$ is a unital ring, we denote by $U^{n\times n}$ the $n\times n$ matrix ring over $U$, and by
$GL(U^{n\times n})$ 
the group of invertible elements of $U^{n\times n}$. If in addition $U$ is a topological ring (e.g. unital
commutative Banach algebra), then we also denote by
$GL_0(U^{n\times n})$
the connected component of $I_n$ (the $n \times n$ identity matrix) in $GL(U^{n\times n})$.

Let $C(G)$ be the algebra of (complex valued) continuous functions on $G$,
and let $P(G)$ be the (non-closed) subalgebra of $C(G)$ of all finite linear combinations of functions $\langle j,
\cdot\rangle$, $j\in\Gamma$,
where
$\langle j,g\rangle$ stands for the action of the character
$j\in\Gamma$ on the group element $g\in G$ (thus, $\langle
j,g\rangle \in \mT$).
 Since $\Gamma$ is
written additively and $G$ multiplicatively, we have
$$ \langle \alpha +\beta ,g\rangle= \langle \alpha,g\rangle \cdot
\langle \beta,g\rangle, \qquad \alpha,\beta\in\Gamma, \ \ g \in G,$$
$$\langle\alpha,gh\rangle=\langle\alpha,g\rangle\cdot\langle\alpha,h\rangle,
\qquad\alpha\in\Gamma,\quad g,h\in G.$$ In the case
$\Gamma={\mathbb Z}$, $G=\mT$ (the unit circle), we will also use
the familiar notation $\langle j, e^{it}\rangle=e^{ijt}$, $0\leq t<2\pi$,
$j\in {\mathbb Z}$. Note that $P(G)$ is dense in $C(G)$ (indeed, by the
Pontryagin duality --- $G$ is the dual of $\Gamma$ --- $\Gamma$
separates points in $G$; now use the Stone-Weierstrass theorem). For
$$a=\sum_{k=1}^m\,a_{j_k}\langle j_k,.\rangle\in P(G), \quad
j_1,\ldots, j_k \in \Gamma \ \ \mbox{are distinct};\, \quad a_{j_k}\neq
0, \ \ k=1,2,\ldots, m,$$    the {\em Bohr-Fourier spectrum} is defined as
the finite set
$$ \sigma (a):=\{j_1, \ldots, j_k\}, $$
and the $a_{j_k}$'s are the  {\em Bohr-Fourier coefficients} of $a$. The
notions of Bohr-Fourier coefficients and Bohr-Fourier spectrum are
extended from functions in $P(G)$ to $C(G)$ by continuity. The
Bohr-Fourier spectrum of $A=[a_{ij}]_{i,j-1}^n\in C(G)^{n\times n}$ is,
by definition, the union of the Bohr-Fourier spectra of the $a_{ij}$'s.
Note that the Bohr-Fourier spectra of elements of $C(G)$ are at most
countable; a proof for the case $\Gamma = \mR$ is found, for
example, in \cite[Theorem 1.15]{C}; it can be easily extended to
general connected compact abelian groups $G$.

We say that a unital Banach algebra ${\cal B}\subseteq C(G)$ is {\em
admissible} if the following properties are satisfied:
\begin{itemize}
\item[(1)]
$P(G)$ is dense in
${\cal B}$;  \item[(2)] ${\cal B}$ is
inverse closed (i.e. $X\in {\cal B}\cap GL(C(G))$ implies $X\in GL({\cal B})$).
\end{itemize}
It is easy to see that property (2) of admissible algebras extends to matrix functions:
\begin{equation}\label{may251}  {\cal B}^{n\times n}\cap GL(C(G)^{n \times n})=GL({\cal B}^{n\times n}).
\end{equation}
Indeed, $A\in  GL(C(G)^{n \times n})$ if and only if $A\in C(G)^{n \times n}$ and
${\rm det}\, A\in GL(C(G))$.

Since $P(G)$ is dense in $C(G)$, it follows that an admissible algebra
is dense in $C(G)$. Important examples of admissible algebras are
$C(G)$ itself and the {\em Wiener algebra} $W(G)$ that consists of all
functions $a$ on $G$ of the form
\begin{equation}\label{nw}
{a}(g)=\sum_{j\in\Gamma}\,a_j\langle j,g\rangle,\qquad g\in G,
\end{equation}
where $a_j\in {\mathbb C}$ and $\sum_{j\in \Gamma}\, |a_j| <\infty$. The norm in $W(G)$ is defined by
$$ \|a\|_1=\sum_{j\in \Gamma}\, |a_j| <\infty. $$
The inverse closed property of $W(G)$ follows from the
Bochner-Philips theorem \cite{BoPh} (a generalization of the classical
Wiener's theorem from the case when  $G=\mT$).

For an admissible algebra ${\cal B}$, we denote by ${\cal B}_\pm$ the
closed unital subalgebra of ${\cal B}$ formed by elements of
${\cal B}$ with the Bohr-Fourier spectrum in $\Gamma_\pm$
(the closedness of ${\cal B}_\pm$ follows from the lower semicontinuity of
the Bohr-Fourier spectrum: if $a_m \rightarrow a$ in ${\cal B}$ as $m\rightarrow \infty$ and
$\sigma (a_m)\subseteq K$ for all $m$, where $K\subseteq \Gamma$ is fixed, then also
$\sigma (a)\subseteq K$).
Thus, ${\cal B}_\pm ={\cal B} \cap C(G)_{\pm}$. Also,
\begin{equation}\label{may196} GL(C(G)_\pm^{n \times n})\cap {\cal B}^{n \times n}=
 GL({\cal B}_\pm^{n \times n}). \end{equation}
Indeed, $\supseteq$ is is obvious. For the converse, if $A\in  GL(C(G)_+^{n \times n})\cap {\cal B}^{n \times n}$, then
$A^{-1}\in C(G)_+^{n \times n}\cap {\cal B}^{n \times n}\subseteq {\cal B}_\pm^{n \times n}$, where
the
inclusion
follows from
(\ref{may251}).

Next, we recall the concept of factorization in the compact abelian
group setting, see e.g. \cite{MRSW,MRS05,EMRS}.  Let ${\cal B}$ be an admissible algebra, and let $A\in
{\cal B}^{n\times n}$. A representation of the form
\begin{equation}\label{rod0}
{A}(g)={A}_-(g)\left({\rm
diag}\,(\langle j_1,g\rangle,\ldots, \langle j_n,g \rangle)\right){A}_+(g),
\quad g\in G,
\end{equation}
where $A_{\pm}, A_{\pm }^{-1}\in {\cal B}_{\pm}^{n \times n}$ and $j_1,
\ldots, j_n\in \Gamma$, is called a (right - with respect to the order $\preceq$) ${\cal B}$-{\em
factorization} of
$A$.
Since $C(G)_+\cap C(G)_-$, and hence also
 ${\cal B}_{+}\cap  {\cal B}_{-}$,  consists of constants only, it follows
that the elements $j_1,\ldots, j_n$ in (\ref{rod0}) are uniquely determined by $A$, up to a permutation. The element
$j_1+\cdots
+j_n\in \Gamma$ is the {\em mean motion} of $A$, denoted ${\rm
MM}\, (A)$. The factorization (\ref{rod0}) is called {\em canonical} if all
$j_1, \ldots, j_n$ are zeros.
Note that (\ref{rod0}) yields (by taking determinants) a ${\cal B}$-factorization of
${\rm det}\, (A(g))$; in particular, $ {\rm
MM}\, (A)= {\rm
MM}\, ({\rm det}\, (A)). $

We say that $A\in{\cal B}^{n\times n}$ is {\em ${\cal B}$-factorable},
resp. {\em canonically $\cal B$-factorable},  if a ${\cal
B}$-factorization, resp. canonical ${\cal B}$-factorization, of $A$ exists.
Clearly, it is necessary that $A\in GL({\cal B}^{n\times n})$ for $A$ to
be ${\cal B}$-factorable.

If $\Gamma=\mZ$, then $G=\mT$,  and $W(\mT)$-factorization is the
classical Wiener-Hopf factorization on the unit circle. As it happens, in
this case the above mentioned necessary invertibility condition is
sufficient as well. This result is due to Gohberg-Krein \cite{GK58}, and
can also be found in many monographs, e.g. \cite{CG,LS}. It is well
known that the condition $A\in GL(C(\mT)^{n\times n})$ is not
sufficient for $C(T)$-factorization even when $n=1$; an
example can be found e.g. in \cite{GF74}.

For $\Gamma=\mR$ the dual group $G$ is the Bohr compactification
of $\mR$, so that $C(G)$ is nothing but the algebra $AP$ of Bohr
almost periodic functions while $W(G)$ is its (non-closed) subalgebra
$APW$ of $AP$ functions with absolutely convergent Bohr-Fourier
series. Characters $\scal{j,.}$ in this setting are simply the
exponentials $e_j(x)=\exp(ij x)$. The $\cal B$-factorization
corresponding to these cases, called $AP$ and $APW$ factorization,
respectively, in the matrix setting was first treated in
\cite{KarlSpit85,KarlSpit86}. It was then observed that there exist
matrix functions in $GL(APW^{2\times 2})$ which are not
$AP$-factorable. In other words, the necessary invertibility condition in
general is not sufficient. The detailed construction of such matrix
functions can be found in \cite{BKS1}, while more recent new classes
are discussed in \cite{CDKS,CaKaS09,RRS09}.

In each of these classes, however, the resulting sets of non $APW$
factorable matrix functions are nowhere dense. This fact has lead
many researchers to view the following conjecture as plausible: {\em
For every admissible algebra ${\cal B}$, the set $GLF ({\cal
B}^{n\times n})$ of ${\cal B}$-factorable elements in ${\cal
B}^{n\times n}$ is dense in $GL({\cal B}^{n\times n})$.} This is the
main issue we are addressing in the present paper.

The conjecture holds for scalar valued functions; this is part of
Theorem \ref{may111} below. (In order not to interrupt the main
stream of the paper, the scalar case is relegated to Section 7.)

Our main,  and rather surprising, finding is that the conjecture is not
true in the matrix case. Namely, the conjecture fails for any $\Gamma$
that contains a copy of ${\mathbb Z}^3$. Moreover, we prove that in
this case the minimal closed subgroup of $GL({\cal B}^{n \times n})$
containing the set of ${\cal B}$-factorable matrix functions, is not
dense in  $GL({\cal B}^{n \times n})$ (Theorem \ref{mar67'}). Even
more, this non-denseness result can be extended to the set of
triangularizable matrix functions (Theorem \ref{mar67''}).
Note also that the set $GLF ({\cal B}^{n\times n})$ is not
necessarily open. For example, it was shown in \cite{CDKS} that
triangular matrix functions \[ \left[\begin{matrix} e_\lambda & f \\ 0 &
e_{-\lambda}\end{matrix}\right] \] with $f\in APW$ such that
$\sigma(f)$ is disjoint with some interval $I\subset
(-\lambda,\lambda)$ of length $d\geq \lambda$ are $APW$ factorable
if and only if the endpoints of $I$ belong to $\sigma(f)$. Consequently,
for $d>\lambda$ such matrix functions lie on the boundary of
$GLF(APW^{2\times 2})$.

Throughout the paper we denote by
${\rm diag}\, (X_1, \ldots , X_m)$ the block diagonal matrix with the diagonal blocks
$X_1,\ldots, X_m$, in that order.

\section{Description of results}
\setcounter{equation}{0}

Our main result is as follows.

\begin{Tm} \label{mar67'}  Let $\Gamma$ be a torsion free abelian group that contains a subgroup
isomorphic to ${\mathbb Z}^3$, and let ${\cal B}$ be an admissible
algebra. Then, for every natural $n\geq 2$ there exist infinitely many
pathwise connected components of $GL({\cal B}^{n \times n})$ with
the property that each one of these components does not intersect
the minimal closed subgroup of $GL({\cal B}^{n \times n})$
containing
 $GLF ({\cal B}^{n \times n})$; in particular, the minimal closed subgroup of $GL({\cal B}^{n \times n})$ containing
 $GLF ({\cal B}^{n \times n})$ is not dense in  $GL({\cal B}^{n \times n})$.
 \end{Tm}

In Theorem \ref{mar67'}, the closedness of a subgroup is understood relative to  $GL({\cal B}^{n
\times n})$.

Of course, the non-denseness of ${\cal B}$-factorable matrix
functions is its easy consequence.

Note that the result of Theorem \ref{mar67'} does not depend on the order $\preceq$.

\bigskip

The proof of Theorem~\ref{mar67'} is contained in Sections~3--6.
Before  passing to their more detailed description, note that in
contrast to Theorem~\ref{mar67'}, for the group of rationals and
its
subgroups, the ${\cal B}$-factorable matrix functions are dense:

\begin{Tm}
If $\Gamma$ is isomorphic to a subgroup of the additive group of
rational numbers ${\mathbb Q}$ (with the discrete topology),
and ${\cal B}$ is an admissible algebra, then $GLF ({\cal B}^{n\times
n})$ is dense in ${\cal B}^{n\times n}$.
\end{Tm}

{\bf Proof.} We can approximate any given element
of ${\cal B}^{n\times n}$ with elements in $P(G)^{n\times n}$ (property (1) of the definition of an
admissible algebra).
Now use the property that
the subgroup of ${\mathbb Q}$ generated by the finite Bohr-Fourier spectrum of any element of $P(G)^{n \times n}$ is
isomorphic to ${\mathbb Z}$.
Thus, it suffices to show that, for the case $\Gamma=\mZ$, $G=\mT$, if
$$ a(t)= \sum_{j=j_0}^{j_1} a_j e^{ijt}, \ \ 0\leq t<2\pi, \ \ a_{j_0}, a_{j_0+1}, \ldots , a_{j_1} \in {\mathbb C}^{n\times
n}, $$
then for some integer $j_2\geq j_1$ the following property holds: for every $\epsilon>0$ there exist
$a'_k\in {\mathbb C}^{n \times n}$, $k=j_0, \ldots, j_2$, such that
$\|a'_k-a_k\|<\epsilon$ for $k=j_0, \ldots, j_2$ (we take $a_{j_1+1}=\cdots =a_{j_2}=0$) and the matrix function
$ a'(t):=\sum_{j=j_0}^{j_2} a'_j e^{ijt} $ is ${\cal B}$-factorable. We may assume without loss of generality
that $j_0=0$. Using the Smith form (the diagonal form) of the matrix polynomial
$$ \widehat{a}(z):=\sum_{j=0}^{j_1} a_j z^j, $$
we see that for some integer $j_2\geq j_1$, for every $\epsilon>0$ there exist
$a'_k\in {\mathbb C}^{n \times n}$, $k=0,1, \ldots, j_2$, such that
$\|a'_k-a_k\|<\epsilon$ for $k=j_0, \ldots, j_2$ and ${\rm det}\, a'(t)\neq 0$ for all $t\in [0,2\pi)$.
The standard procedure for factorization of rational matrix functions (see e.g. \cite[Section XIII.2]{GGK1}) now implies that
$a'$ is ${\cal B}$-factorable.
\  \ \ \ \ \ $\Box$
\medskip

For groups $\Gamma$ that are not subgroups of ${\mathbb Q}$ and at the same time do not contain
${\mathbb Z}^3$ (example: $\Gamma={\mathbb Z}^2$), we do not know whether or not $GLF ({\cal B}^{n\times n})$ is dense in
$GL({\cal B}^{n\times n})$,
where ${\cal B}$ is an admissible algebra. In this respect note that every pathwise
connected component
of $GL(C({\mathbb T}^2)^{n \times n})$ contains
$C({\mathbb T}^2)$-factorable elements of $C({\mathbb T}^2)^{n \times n}$, in contrast to Theorem
\ref{mar41}. This can be easily seen by comparing Proposition \ref{dec183} and Theorem \ref{sep181}; indeed,
the  pathwise
connected components of $GL(C({\mathbb T}^2)^{n \times n})$, as well as those of the set of
$C({\mathbb T}^2)$-factorable elements of $C({\mathbb T}^2)^{n \times n}$, are parametrized identically by
${\mathbb Z}^2$.

In the next section we develop some preliminary results on  pathwise
connected components that are perhaps of independent interest. The
proof of a preliminary result and some corollaries are given in Section
5. In Section 4 we present an example treating in detail the key case of
the three-dimensional torus. In Section 6 we study the minimal
subgroup containing all factorable matrix functions and complete the
proof of Theorem \ref{mar67'}. As was mentioned earlier, the
scalar case is studied in Section 7. Finally,  in Section 8 we extend our
main result to non-density of triangularizable matrix functions.

\section{Connected components: finite dimensional tori}
\setcounter{equation}{0}

Consider the following problem:
\begin{Pb}\label{dec181}
Describe the  pathwise connected components of the group
$GL({\cal B}^{n \times n})$, where ${\cal B}$ is an admissible algebra.
\end{Pb}

Note that the set of pathwise connected components of
$GL({\cal B}^{n \times n})$ has a natural group structure: If
${\cal C}_1, {\cal C}_2$ are two such components, then
${\cal C}_1 \cdot {\cal C}_2$ is the component defined by the property that
$A_1A_2\in {\cal C}_1 \cdot {\cal C}_2$, where $A_1\in {\cal C}_1$, $A_2\in {\cal C}_2$.

To study Problem \ref{dec181}, the following proposition will be handy.
\begin{Pn} \label{mar61} Let ${\cal B}$ be an admissible algebra. Then:

${\rm (1)}$  The pathwise connected components of
$GL({\cal B}^{n \times n})$ have the form
\begin{equation}\label{mar62} GL({\cal B}^{n \times n})\cap {\cal C}, \qquad {\cal C}\in C_n(G), \end{equation}
where $C_n(G)$ stands for the group of pathwise connected components of $GL(C(G)^{n
\times n})$.

${\rm (2)}$ $GL({\cal B}_{\pm}^{n \times n})$ is pathwise connected.
\end{Pn}

{\bf Proof.} Part (1) follows from Arens' theorem \cite{Ar66}; we give an independent
selfcontained proof. Clearly, if
$f_1:[0,1] \rightarrow GL({\cal B}^{n \times
n})$ is a continuous path, then
$f_1(0)$ and $f_1(1)$ belong to the same connected component of $GL(C(G)^{n \times n})$. Thus, any pathwise connected component
of  $GL({\cal B}^{n \times n})$ is contained in one of the sets (\ref{mar62}). Next, we
prove that
if $X,Y\in  GL({\cal B}^{n \times n})\cap {\cal C}$ for some ${\cal C}\in C_n(G)$, then $X$ and $Y$ are pathwise connected in
$GL({\cal B}^{n \times n})$. Let $f_2: [0,1] \rightarrow GL({C(G)}^{n \times n})$ be a continuous path such that
$f_2(0)=X$, $f_2(1)=Y$. We necessarily have $f_2([0,1])\subset {\cal C}$. Since $GL({C(G)}^{n \times n})$ is open, there exists
$\va>0$ such that
\eq{mar66} Z\in {C(G)}^{n \times n}, \ \ \|Z-f_2(t)\|_\infty<\va \ \ \mbox{for some $t\in [0,1]$}
\quad \Longrightarrow \quad Z\in GL({C(G)}^{n \times n}). \en
Let $m$ be a sufficiently large integer so that
\begin{equation}\label{mar63} \|f_2((j+1)/m)-f_2(j/m)\|_\infty < \va, \qquad j=0,1,\ldots, m-1.
\end{equation}
Let
\eq{mar64} \kappa:=\max_{j=0,1,\ldots, m-1} \{\|f_2((j+1)/m)-f_2(j/m)\|_\infty\}<\va, \en
and select $\lambda>0$ such that $\kappa+\lambda<\va$. Since ${\cal B}$
is admissible, we can select $F_{j/m}\in {\cal B}^{n \times n}$, $j=0,1,\ldots, m$ such that
$F_0=X$, $F_1=Y$, and $\|F_{j/m}-f_2(j/m)\|_\infty <\lambda$ for $j=1,2,\ldots, m-1$.
Then (\ref{mar63}), (\ref{mar64}) guarantee that
$$ [F_{j/m}, F_{(j+1)/m}]:=
\{tF_{j/m} +(1-t)F_{(j+1)/m}\, : \, 0\leq t\leq 1\}$$ is contained in $$
\{Z\in C(G)\, :\, \|Z-f_2((j+1)/m)\|_\infty <\va\} \cap \{Z\in C(G)\, :\, \|Z-f_2(j/m)\|_\infty <\va\},
$$
for $j=0,1,\ldots, m-1$. By (\ref{mar66}), the interval $ [F_{j/m}, F_{(j+1)/m}]$ is contained in
$GL({C(G)}^{n \times n})$, and since ${\cal B}$ is admissible, we actually have
 $ [F_{j/m}, F_{(j+1)/m}]\subset  GL({\cal B}^{n \times n})$. Now clearly
the union $\cup_{j=0}^{m-1}[F_{j/m}, F_{(j+1)/m}]$ represents a continuous path in
$ GL({\cal B}^{n \times n})$  connecting $X$ and $Y$.

Finally, observe that $GL({\cal B}^{n \times n})$ is dense in
$GL(C(G)^{n\times n})$ (as easily follows from the defining properties of
admissible algebras), and therefore $GL({\cal B}^{n \times n})\cap {\cal C}\neq \emptyset$
for every ${\cal C}\in C_n(G)$.
\medskip

Part (2). As in Part (1), using (\ref{may196}) we show that the
pathwise connected components of
 $GL({\cal B}_{\pm}^{n \times n})$ are intersections of $GL({\cal B}_{\pm}^{n \times n})$ with the
pathwise connected components of $GL(C(G)_{\pm}^{n \times n})$.
But the latter is pathwise connected \cite[Theorem 6.1]{BrudRS},
and so we are done.  $\ \ \ \ \Box$
\bigskip

In view of Proposition \ref{mar61}, we may assume ${\cal B}=C(G)$ in Problem \ref{dec181}.
In the classical case $\Gamma=\mZ$ the answer is known: The group
$GL(C(G)^{n \times n})$ has infinitely many  pathwise connected components each
of which is characterized by a fixed integer -- the winding number of
the determinant (see, e.g., Bojarski's appendix to \cite{Mus68}).
In general, $C_n:=GL(C(G)^{n \times n})/GL_0(C(G)^{n \times n})$
can be identified with
the group $[G;GL(\mC^{n\times n})]$ of homotopy classes of continuous maps $G\to GL(\mC^{n\times n})$.

In this section, we give a solution of Problem \ref{dec181} for the case when $G$ is a finite dimensional torus. The solution
for
an infinite dimensional $G$ is obtained from here using the fact that $G$ is the inverse limit of a
family $\{G_\alpha\}_{\alpha\in\Lambda}$ of finite dimensional tori and hence the set of connected
components of $GL(C(G))$ in this case is the direct limit of sets of
connected components of $GL(C(G_\alpha))$ under the maps transposed to the corresponding maps of the inverse
limit.
The solution for $G$ being a finite dimensional torus is based on the classical results of Fox \cite{Fox} concerning torus homotopy groups and their
connections with certain homotopy groups of the unitary group $U_n\subset GL(\mC^{n\times n})$.

For the readers' convenience, we recall basic definitions and some results from \cite{Fox}. We identify the $k$-dimensional
torus $\mT^k$ with $\mR^k/\mZ^k$. Let $E=[0,1]$ and let $E^k$ stand for the unit $k$-dimensional cube.
Let $X$ be a topological space. We say that a continuous function $f:E^k\rightarrow X$ has {\em periodic boundary conditions}
if for each $i=1,2,\ldots, k$,
$$ f(\ldots,x_{i-1},0, x_{i+1}, \ldots )=  f(\ldots,x_{i-1},1, x_{i+1}, \ldots ) \qquad \forall
\ \ 0\leq x_j\leq 1, \ (j\neq i). $$
Clearly, there is a one-to-one correspondence between such functions $f$ and continuous functions
$F: \mT^k \rightarrow X$.

Let $o\in X$ be a fixed point. We denote by ${\cal T}^k(X,o)$ the set of all continuous functions
 $f:E^k\rightarrow X$ with periodic boundary conditions which satisfy in addition the condition
$$ f(0,x_2\ldots, x_k)=f(1,x_2\ldots, x_k)=o. $$
By $\tau_k(X,o)$ we denote the set of all homotopy classes of functions
$f\in {\cal T}^k(X,o)$. Then the {\em torus homotopy group} is the set $\tau_k(X,o)$ together
with the group binary operation ``+''
induced by the operation of multiplication
$$ (f_1\times f_2)(x)=\left\{\begin{array}{ll} f_1(2x_1,x_2,\ldots, x_k) & \mbox{if $0\leq x_1\leq 1/2$,} \\[3mm]
f_2(2x_1-1,x_2,\ldots, x_k) & \mbox{if $1/2\leq x_1\leq 1$,} \end{array}\right.
\qquad (f_1,f_2\in  {\cal T}^k(X,o)), $$
with the inverse operation being
$$ f^{[-1]}(x)=f(1-x_1, x_2,\ldots, x_k), \qquad f\in {\cal T}^k(X,o). $$
If $X$ is pathwise connected, then  $\tau_k(X,o)$ does not depend on
the particular choice of $o$ (up to a group isomorphism).

\begin{Pn} {\rm \cite{Fox}} If $X$ is a topological group with the unit element $e$, then
$(\tau_k(X,e),+)$
is an abelian group.
 \end{Pn}

In the sequel we will work with the torus homotopy groups of $U_n$
and use the notation $\tau_k(U_n)$ for
$(\tau_k(U_n,I),+)$. One can introduce another multiplication
``$\,\cdot\,$'' on $\tau_k(U_n)$ induced from the matrix multiplication
in $U_n$.  Note that $\tau_1(U_n)$ is the fundamental group of $U_n$,
and so $\tau_1(U_n)\cong\mZ$. In this case one also has
$\tau_1(U_n)\cong (\tau_1(U_n),\cdot)$ but in general for $k,n\ge 2$
these two groups are not isomorphic. In what follows we will use the
notation $\dot\tau_k(U_n):=(\tau_k(U_n),\cdot)$. Also, we will use the well-known semidirect products of groups.
Briefly, a group $X$ is a
semidirect product of its normal subgroup $N$ and a subgroup $Y$,
notation: $X=Y \ltimes N$, if $X=NY=YN$ and the intersection of $N$
and $Y$ is trivial; see e.g. \cite{DF} for more information.

\begin{Pn}\label{dec183}  Assume $\Gamma=\mZ^k$, i.e. $G=\mT^k$, and let ${\cal B}$ be an admissible algebra. Then
the group of  pathwise connected components of $GL({\cal B}^{n\times n})$, $n\geq 1$,
equipped with the multiplication induced by the product of matrices in $GL_n(\mathbb C)$, is isomorphic
to
$\dot\tau_k(U_n)\ltimes(\dot\tau_{k-1}(U_n)\ltimes\cdots\ltimes(\dot\tau_2(U_n)\ltimes\dot\tau_1(U_n))\dots)$
with naturally defined multiplications in the semidirect products. If $k\le 2n-1$, then the direct product of
groups $\tau_k(U_n)\times \tau_{k-1}(U_n)\times\cdots\times\tau_1(U_n))$ is isomorphic to $\mZ^{2^{k-1}}$.
\end{Pn}

{\bf Proof.} In view of Proposition \ref{mar61}, we may assume ${\cal B}=C(G)$. Note that $U_n$ is
a strong deformation retract
\cite{Hust1} of $GL(\mC^{n
\times n})$;
indeed, the map $$(X,t) \ \mapsto \ X(X^*X)^{-1/2}(t(X^*X)^{1/2}+(1-t)I), \quad X\in GL(\mC^{n \times n}),
\ \
0\leq
t\leq 1, $$ provides such a retracting deformation. Therefore, the sets of
 pathwise connected components
of $GL(C(\mT^k)^{n\times n})$ and
of $U(C(\mT^k)^{n\times n})$ ($:=$ the set of continuous functions on $\mT^k$ with values in $U_n$) coincide.

By $x_1,\dots, x_k$ we denote the local coordinates on $\mT^k$
induced from the coordinates on $\mR^k$ by means of the quotient
map $\mR^k\to\mR^k/\mZ^k=\mT^k$. Consider the subtorus
$\mT^{k-1}$ of codimension $1$ in $\mT^k$ defined by
$\{(x_1,x_2,\dots, x_k)\in\mT^k\, :\, x_k=0\}$. By
$p_{k-1}:\mT^k\to\mT^{k-1}$ we denote the projection sending a
point with coordinates $(x_1,\dots, x_k)$ to the point with coordinates
$(x_1,\dots, x_{k-1}, 0)$. Each function $X\in U(C(\mT^k)^{n\times n})$
can be uniquely written as
\[
X=X_1\cdot X_2,\quad\text{where}\quad X_2:=(X|_{\mT^{k-1}})\circ p_{k-1},\quad X_1:=X\cdot X_2^{-1}.
\]
The matrix $X_1$ belongs to the Banach group
$U_I(C_1(\mT^k)^{n\times n})$ of continuous functions on $\mT^k$
with values in $U_n$ identically equal  $I_n$ on $\mT^{k-1}$. Here
$C_1(\mT^k)$ is the Banach algebra of (complex) continuous
functions on $G:=\mT^k$ which are constant on $\mT^{k-1}$. (Recall
that a Banach group is a Lie group that is modelled locally by open
balls in a Banach space.) Clearly, $U_I(C_1(\mT^k)^{n\times n})$ is a
normal subgroup of $U(C(\mT^k)^{n\times n})$ with the quotient
group $U(C(\mT^{k-1})^{n\times n})$. Moreover the above
decomposition shows that the exact sequence of Banach groups
\begin{equation}\label{eq7.1}
1\longrightarrow U_I(C_1(\mT^k)^{n\times n})\longrightarrow U(C(\mT^k)^{n\times n})\longrightarrow U(C(\mT^{k-1})^{n\times n})\longrightarrow 1
\end{equation}
splits, i.e. $U(C(\mT^k)^{n\times n})$ is a semidirect product of subgroups $U_I(C_1(\mT^k)^{n\times n})$ and
$$(p_{k-1})^*U(C(\mT^{k-1})^{n\times n})\cong
U(C(\mT^{k-1})^{n\times n}). $$
Here, $p_{k-1}^*\, :\, U(C(\mT^{k-1})^{n\times n}) \ \rightarrow \ U(C(\mT^{k})^{n\times n})$
 is defined by
$$ (p_{k-1}^*(f))(x)=f(p_{k-1}(x)), \quad f\in  U(C(\mT^{k-1})^{n\times n}), \quad        x\in \mT^k. $$
Recall that the set of  pathwise connected components of a Banach group $Q$ can be
naturally
identified with the (discrete) group defined as the quotient of $Q$ by the  pathwise connected component containing the
identity
(which is a normal subgroup of $Q$). In the case of $U(C(\mT^k)^{n\times n})$, the group of
 pathwise connected components $C_{k,n}$ is a discrete countable group (because the former group is separable).
Denote by $\widetilde C_{k,n}$ the normal subgroup of $C_{k,n}$ isomorphic to the group of  pathwise
connected components of
the
Banach group
$U_I(C_1(\mT^k)^{n\times n})$.
Note that we have a split exact
sequence
$$ 1 \ \longrightarrow \ \widetilde C_{k,n} \ \longrightarrow \ C_{k,n} \ \longrightarrow
\ C_{k-1,n} \ \longrightarrow \ 1. $$
This follows from  \eqref{eq7.1} and the fact that the connected component of $U(C(\mT^k)^{n\times n})$ containing the
unit matrix is the semidirect product (under the multiplication given by product of matrices) of the connected components
containing the unit matrix of groups  $U(C(\mT^{k-1})^{n\times n})$ and $U_I(C_1(\mT^k)^{n\times n})$. Moreover, this
connected component of  $U_I(C_1(\mT^k)^{n\times n})$ is a normal subgroup of $U(C(\mT^k)^{n\times n})$.
Thus, $C_{k,n}$ is isomorphic to a semidirect product of discrete
subgroups $\widetilde C_{k,n}$ and $C_{k-1,n}$, where
$\widetilde C_{k,n}$ is a normal subgroup of $C_{k,n}$ isomorphic to the group of  pathwise connected components of the
Banach group
$U_I(C_1(\mT^k)^{n\times n})$, and $C_{k-1,n}$ is isomorphic to the group of  pathwise connected components of
$U(C(\mT^{k-1})^{n\times n})$.
We parametrize elements of $C_{k,n}$ by points of $\widetilde C_{k,n}\times C_{k-1,n}$.
Namely, each element of $C_{k,n}$ can
be uniquely
presented as a product of elements of $\widetilde C_{k,n}$ and $C_{k-1,n}$ (taken in this order) and so $C_{k.n}$ can
be
recovered from
$\widetilde C_{k,n}\times C_{k-1,n}$.

Next,
the set of  pathwise connected components of $U(C_1(\mT^k)^{n\times n})$ is
naturally identified with the set of
homotopy classes of continuous maps of the maximal ideal space $M_1$ of the Banach algebra $C_1(\mT^k)$ into $U_n$. By the definition $M_1$
is homeomorphic to the compact Hausdorff space obtained from $\mT^k$ by contracting the subtorus $\mT^{k-1}$ into a point. Therefore the set of
 pathwise connected components of $U_I(C_1(\mT^k)^{n\times n})$ is identified with the homotopy classes of continuous
maps of $\mT^k$ into $U_n$ which map
$\mT^{k-1}$ into $I_n\in U_n$. Then, by the definition of the multiplication in $\widetilde C_{k,n}$ we have
$\widetilde C_{k,n}\cong\dot\tau_k(U_n)$.

Applying the same arguments to $U(C(\mT^{k-1})^{n\times n}),\dots, U(C(\mT^{1})^{n\times n})$ we obtain the first statement of the proposition. In particular,
each  pathwise connected component of \penalty-10000 $GL(C(\mT^k)^{n\times n})$ has the form $X_k\cdots X_1\cdot
GL_0(C(\mT^k)^{n\times n})$,
where $X_r$ represents an element of ${\cal T}^r(U_n,I)$, $1\le r\le k$, and
$GL_0(C(\mT^k)^{n\times n})$ is the
pathwise connected component containing $I_n$.

Further, it is known that $\tau_r(U_n)$, $r\ge 2$, is
defined from the split exact sequence of abelian groups, see \cite[Statement (9.3)]{Fox},
\begin{equation}\label{e7.2}
0\longrightarrow \prod_{i=2}^r (\pi_i(U_n))^{\alpha_i}\longrightarrow \tau_r(U_n)\longrightarrow\tau_{r-1}(U_n)\longrightarrow 0,
\end{equation}
where $\alpha_i:={r-2 \choose i-2}$, $\pi_i(U_n)$ is the $i$-dimensional homotopy group of $U_n$ and
$(\pi_i(U_n))^{\alpha_i}$ is the $\alpha_i$-fold direct product $\pi_i(U_n)\times\cdots\times\pi_i(U_n)$.

To compute $\pi_i(U_n)$ one can use e.g. the Bott periodicity theorem \cite{Bott} which, in particular, implies that
\begin{equation}\label{e7.3}
\pi_{i}(U_n)=\pi_{i+2}(U_{n+1})\qquad\text{for}\qquad 0\leq i\le 2n-1.
\end{equation}
For instance, using the fact that $\pi_0(U_n)=0$ (because $U_n$ is  pathwise connected) and  $\pi_2(U_n)=0$
(see e.g. \cite[Section 8.12]{Huse}),
we obtain from the latter
periodicity that the set of  pathwise connected components of $GL(C(\mT^1)^{n\times n})$ can be naturally
identified with
$\tau_1(U_n)\cong\mathbb Z$ (as remarked above), the set of
 pathwise connected components of $GL(C(\mT^2)^{n\times n})$ can be naturally identified with
$\tau_2(U_n)\times\tau_1(U_n)\cong \mathbb Z\times\mathbb Z$ \cite{Wong}, for $n\ge 2$ the
set of  pathwise connected components of $GL(C(\mT^3)^{n\times n})$ can be naturally identified with
$\tau_3(U_n)\times\tau_2(U_n)\times\tau_1(U_n)\cong (\mathbb Z\times\mathbb Z)\times\mathbb Z\times\mathbb
Z=\mathbb Z^4$,  whereas if $n=1$ it is identified with $\mathbb Z^3$, etc.

In particular, if $k\le 2n-1$, then the set of  pathwise connected components of $GL(C(\mT^k)^{n\times n})$ is
parametrized by
$\tau_k(U_n)\times\cdots\times\tau_1(U_n)\cong\mZ^{2^{k-1}}$.
This easily follows by induction from \eqref{e7.2}, \eqref{e7.3} using the fact that the sum of all odd binomial coefficients
among ${r-2 \choose i-2}$ is equal to $2^{r-3}$, $r\ge 3$.\qquad $\Box$
\bigskip

\section{Example}
\setcounter{equation}{0}

In the case of $\mathbb T^3$ we give an explicit description of the
identification in Proposition \ref{dec183}. This will be later used in the
proof of Theorem \ref{mar41}.
\begin{Ex}\label{e2.4}
{\rm As in the proof of Proposition \ref{dec183}, we
denote by $U_I(C_1(\mT^3)^{n\times n})$ the Banach group of
continuous functions on $\mT^3$ with values in $U_n$ assuming the
value $I_n$ on $\mT^{2}$.
 Let $X\in U_I(C_1(\mathbb
T^3)^{n\times n})$, $n\ge 2$, be a
continuous $n \times n$ unitary matrix function sending
$T_3:=\mathbb
T^2=\{(x_1,x_2,x_3)\in\mathbb T^3\, :\, x_3=0\}\subset\mathbb T^3$ into $I_n\in U_n$. We will consider also subtori
$$T_2:=\{(x_1,x_2,x_3)\in\mathbb T^3\, :\, x_2=0\}\quad \mbox{ and}\quad
T_1:=\{(x_1,x_2,x_3)\in\mathbb T^3\, :\, x_1=0\}$$ of $\mathbb T^3$ with projections
$$t_2:\mathbb T^3\to T_2, \ \ (x_1,x_2,x_3)\mapsto (x_1,0,x_3) \quad \mbox{ and} \quad t_1:\mathbb T^3\to
T_1,
 \ \ (x_1,x_2,x_3)\mapsto
(0,x_2,x_3), $$ respectively. We write
$$ X=X_1\cdot X_2\cdot X_3, \qquad X_j\in  U_I(C_1(\mathbb T^3)^{n\times n}), \ \ j=1,2,3,$$
where  
\begin{eqnarray*} X_3(x_1,x_2,x_3)&=& X(x_1,0,x_3), \\
X_2(x_1,x_2,x_3)&=&X(0,x_2,x_3)X(0,0,x_3)^{-1}, \\
 X_1(x_1,x_2,x_3)&=&X(x_1,x_2,x_3)X(x_1,0,x_3)^{-1} X(0,0,x_3)X(0,x_2,x_3)^{-1}. \end{eqnarray*}
Considering $X_1,X_2,X_3$ as $U_n$-valued periodic functions on $\mathbb R^3$ we obtain that $X_1$ sends the boundary
$\partial E^3$ of the three dimensional closed unit cube $E^3$ into $I_n$; $X_2$ is the pullback by
$t_1$ of some periodic unitary matrix function $\widetilde X_2$ defined on the plane
$R_1:=\{(0,x_2,x_3)\in\mathbb R^3\}$ such that $$\widetilde X_2(0,x_2,0)=\widetilde X_2(0,x_2,1)=\widetilde
X_2(0,0,x_3)=\widetilde X_2(0, 1, x_3)=I_n, $$
i.e.,
$$ X_2(x)=\widetilde X_2(t_1(x)), \qquad x\in {\mathbb T}^3; $$
and $X_3$ is the pullback by $t_2$ of some periodic unitary matrix function
$\widetilde X_3$ defined on the plane $R_2:=\{(x_1,0,x_3)\in\mathbb R^3\}$ such that $\widetilde X_3(x_1,0,0)=\widetilde X_3(x_1,0,1)=I_n$.

Let $\dot{E^3}$ be the interior of the three dimensional closed cube $E^3$. The one-point compactification of $\dot{E^3}$ is homeomorphic to
the unit three dimensional sphere $\mathbb S^3\subset \mathbb R^4$. Moreover, there exists a continuous surjection $E^3\to\mathbb S^3$
which maps the boundary of $\partial E^3$ of the cube into the point $P=(1,0,0,0)$ and
$E^{3}\setminus\partial E^3$ bijectively onto $\mathbb S^3\setminus P$. This map can be extended by periodicity to $\mathbb R^3$ and this extension determines a continuous surjective map
$\psi:\mathbb T^3\to\mathbb S^3$ such that the induced map of the \v{C}ech cohomology groups
$$\psi^* : H^3(\mathbb S^3,\mathbb Z)\cong \mathbb Z \, \to\,  H^3(\mathbb T^3,\mathbb Z)\cong \mathbb Z $$ is
an isomorphism. (That $\psi^*$ is an isomorphism follows from $\psi$ being a bi-continuous bijection
between open dense subsets
$\dot{E^3}\subset\mathbb T^3$ and $\mathbb S^3\setminus P\subset\mathbb S^3$.) By the definition of
the function $X_1$, there exists a
continuous function $X_1'$ on $\mathbb S^3$ with values in $U_n$ whose
pullback by $\psi$ coincides with $X_1$ and such that $X_1'(P)=I_n$. In particular, $X_1'$ determines
an element of the homotopy group $\pi_3(U_n)$. Similarly $\widetilde X_2$ maps the boundary of the
square $E^3\cap R_1$ into $I_n$ and therefore is
the pullback of a function of the two-dimensional unit sphere $\mathbb S^2$ into $U_n$ and so determines an element
of
$\pi_2(U_n)$. But the latter group is trivial and therefore $\widetilde X_2$ considered as a map of the two-dimensional torus into $U_n$ is homotopic to a constant map implying that $X_2$ belongs to the connected component of $U_I(C_1(\mathbb T^3)^{n\times n})$ containing the constant map.

Further, analogously to $X$ the function $\widetilde X_3$ can be factorized as
\[
\widetilde X_3=\widetilde X_{31}\cdot \widetilde X_{32},
\]
where $\widetilde X_{31}$ maps the boundary of the square $E^3\cap R_2$ into $I_n$ and therefore as above (since $\pi_2(U_n)=0$) $\widetilde X_{31}\circ t_2$ belongs to the connected component of $U_I(C_1(\mathbb T^3)^{n\times n})$
containing the constant map, and $\widetilde X_{32}$ is the pullback by $t_1$ of a continuous periodic unitary function $\widetilde X_{32}'$ defined on the line $\{(0,0,x_3)\subset\mathbb R^3\}$. This function by definition determines an element of $\pi_1(U_n)\cong\mathbb Z$. Moreover, each such $\widetilde X_{32}'$ is homotopic to a function
${\rm diag}(1,\dots, 1, e^{2\pi i kx_3})$ for some $k\in\mathbb Z$. Combining the above facts we obtain that
\begin{equation}\label{feb11}
X=X'\cdot {\rm diag}(1,\dots, 1, e^{2\pi i kx_3})\cdot X'',
\end{equation}
where $X'$ is the pullback by $\psi$ of a continuous map $\widetilde X':\mathbb S^3\to U_n$ representing an
element of $\pi_3(U_n)$ and $X''$ belongs
to the connected component of $U_I(C_1(\mathbb T^3)^{n\times n})$ containing the constant map.
Note that $k$ and the element of $\pi_3(U_n)$ in (\ref{feb11}) are
unique, for a given $X$, but the factorization \eqref{feb11} is generally not unique (this fact will not be used in the
sequel).

Let us describe some canonical matrix functions representing the first terms of the above factorization.

Consider the group action of $U_n$ on  $\mathbb S^{2n-1}$ defined by the map
$p:U_n \times \mathbb S^{2n-1} \ \rightarrow\ \mathbb S^{2n-1}$,
$$ p(X,(s_1,\ldots ,s_{2n}))=X\left[\begin{array}{c} s_1+is_2 \\ s_3+is_4 \\ \vdots \\ s_{2n-1}+is_{2n}
\end{array}\right], $$ where $$ X\in U_n, \quad  s_1,\ldots, s_{2n}\in {\mathbb R}, \quad  s_1^2+\cdots+ s_{2n}^2=1. $$
For $e_1:=(1,0,\dots, 0)\in\mathbb S^{2n-1}$ the subgroup of $U_n$
fixing $e_1$ is $1\times U_{n-1}\cong U_{n-1}$. The group $1\times
U_{n-1}$ acts on $U_n$ by left multiplications and the set of
equivalence classes $U_n/(1\times U_{n-1})$ under this action is a
homogeneous space diffeomorphic to $\mathbb S^{2n-1}$. The latter
diffeomorphism  is defined by identifying the first column of each
matrix from $U_n$ with an element of $\mathbb S^{2n-1}$. Now the
composition $\pi$ of maps $U_n\to U_n/(1\times U_{n-1})\to\mathbb
S^{2n-1}$ defines $U_n$ as a principal fibre bundle over $\mathbb
S^{2n-1}$ with fibre $U_{n-1}\cong 1\times U_{n-1}\subset U_n$. In
particular, over any proper open subset $S$ of $\mathbb S^{2n-1}$
there exists a smooth right inverse to the map $\pi$
 (because $\mathbb S^{2n-1}\setminus\{s\}$ is contractible
for each $s\in\mathbb S^{2n-1}$). This implies that
$\pi^{-1}(S)$ is diffeomorphic to $S\times U_{n-1}$. Using this
construction one can easily show (comparing dimensions of
$\mathbb S^3$ and of $\mathbb S^{2n-1}$) that each continuous
map $\mathbb S^3\to U_n$ is homotopic to a continuous map
$\mathbb S^3\to U_2$, where $U_2$ is identified with the subgroup
$1\times 1\times\cdots\times 1\times U_2\subset U_n$.

Moreover, if two continuous maps $f_0,f_1: \mathbb S^3 \to 1 \times U_{n-1}
$ ($n\geq 3$) are homotopic to each other via a homotopy $f_t: \mathbb S^3 \to U_n$, $0 \leq t \leq 1$, then they are also
homotopic to
each other via a homotopy $\widetilde{f}_t: \mathbb S^3 \to 1 \times U_{n-1}$.
Indeed, re-write $f_t$ as a map
$f:[0,1]\times \mathbb S^3 \to U_n$, where
$$ f_0(\mathbb S^3)=f(0, \mathbb S^3)\subseteq 1 \times U_{n-1}, \quad
f_1(\mathbb S^3)=f(1, \mathbb S^3)\subseteq 1 \times U_{n-1}. $$
Using a dimension argument ($n\geq 3$) we may assume that
$\pi \circ f: [0,1]\times \mathbb S^3 \to \mathbb S^{2n -1}$ is not surjective; here $\pi$ is defined in the
preceding paragraph. Using the principal fiber bundle construction as
described there, we see that there is a homotopy
$F:[0,1]\times [0,1] \times \mathbb S^3$ connecting $f$ to some $\widehat{f}: [0,1] \times \mathbb S^3
\to 1 \times U_{n-1}$, i.e. $F(0,\cdot, \cdot)=f$, $F(1,\cdot, \cdot)=\widehat{f}$.
Moreover, from the fibre bundle argument it follows that one can make sure that if
$f(t,z)\in 1 \times U_{n-1}$, then $F(s,t,z)=f(t,z)$ for all $s\in [0,1]$.
Now $\widehat{f}$ is the desired homotopy between $f_0$ and $f_1$.

It follows that the natural map $\tau: \pi_3(U_{n-1}) \to \pi_3(U_n)$ is an isomorphism, and by induction
we see that $\pi_3(U_{2})$ is naturally isomorphic to  $\pi_3(U_n)$.
In turn,
$U_2$ is diffeomorphic to $\mathbb S^1\ltimes \mathbb S^3$ ($U_2$
is a semidirect product of $U_1$ and of the $2 \times 2$ special
unitary group $SU_2\cong \mathbb S^3$). Since $\pi_3(U_2)=
\pi_3(\mathbb S^1\ltimes \mathbb S^3)$ is isomorphic to
$\pi_3(SU_2)=\pi_3(\mathbb S^3)$ (because
$\pi_3(\mathbb S^1)$ is trivial), we obtain
therefore that each
continuous map $\mathbb S^3\to U_n$ is homotopic to a continuous
map $\mathbb S^3\to \mathbb S^3$ under the identification
$\mathbb S^3= SU_2$.

Next, consider the map $\mathbb R^3\to\mathbb R^3$, $(x_1,x_2,x_3)\mapsto (kx_1, x_2, x_3)$,
$k\in\mathbb Z$.
For $k\ne 0$, it determines a finite map $d_k:\mathbb T^3\to\mathbb T^3$; the number of preimages of each point under $d_k$ is $k$ and $d_k$
preserves orientation for $k>0$ and reverses it for $k<0$.
For $k=0$ it determines a projection $d_0:\mT^3\to T_1$, where $T_1:=\{x=(0,x_2,x_3)\in\mT^3\}$
is a two-dimensional subtorus.  We set
\[
\psi_k:=\psi\circ d_k .
\]
 By definitions of $\psi$ and $d_k$ there exists a continuous map $e_k:\mathbb S^3\to\mathbb S^3$ such that
$\psi_k=e_k\circ\psi$, for $k\ne 0$ preimage of each
point of $\mathbb S^3\setminus\{(1,0,0,0)\}$ under
$e_k$ consists of $k$ distinct points, and preimage of
$(1,0,0,0)$ has dimension $\le 2$. Also,
$\psi_0$ has two-dimensional image. By the definition $e_k$
determines a map in cohomology $e_k^*:H^3(\mathbb S^3,\mathbb
Z)\cong\mZ\to H^3(\mathbb S^3,\mathbb Z)\cong\mZ$ such that
$e_k^*(m):=km$, $m\in\mathbb Z$.  Assume now that $\widetilde X'$
as above represents a nonzero element $s\in
\pi_3(U_n)=\pi_3(SU_2)\, (=\pi_3(\mathbb S^3))$. If $k\, (\ne 0)$ is the
degree of $s$, then by the Hopf theorem (as stated for example in
\cite[Section II.8]{Hust1}), $\widetilde X'$ is homotopy equivalent to
$e_k$. From here identifying $\mathbb S^3$ with $SU_2$ by sending
$x=(x_1,x_2,x_3,x_4)\in\mathbb S^3$ into the matrix
$$
S(x):=\left[
\begin{array}{cc}
x_1+ix_2&-(x_3-ix_4)\\
x_3+ix_4&x_1-ix_2
\end{array}
\right]\in SU_2
$$
we obtain that each $X\in U_I(C_1(\mathbb T^3)^{n\times n})$, $n\ge 2$, admits a factorization
\[
X={\rm diag}(1,\dots, 1,(S\circ\psi_m))\cdot {\rm diag}(1,\dots, 1, e^{2\pi i p x_3})\cdot Y
\]
where $m,p\in\mZ$ and $Y$ belongs to the connected component of $U_I(C_1(\mathbb T^3)^{n\times n})$
containing the
constant map. The numbers $m,p$ here are uniquely defined but the above factorization is not unique.
}
\end{Ex}

\section{Connected components without factorable elements}
\setcounter{equation}{0}

In this section we prove the following result which is a particular case of, and
at the same time a
stepping stone to, the
proof of Theorem \ref{mar67'}.

\begin{Tm}\label{mar41} Let $\Gamma$ be
 a torsion free abelian group that contains a subgroup isomorphic to $\mZ^3$.
Then for every natural $n\geq 2$, there exist infinitely many pathwise connected components
of $GL(C(G)^{n \times n})$ with the property that each one of these components does not contain
any $C(G)$-factorable element of $C(G)^{n \times n}$.
\end{Tm}

The following result obtained in \cite{BrudRS} will be used in the proof:

\begin{Tm} \label{sep181} The pathwise connected components of $GLF
(C(G)^{n\times n})$ are precisely the sets
\begin{equation} \label{feb61}  \{X\in GLF (C(G)^{n\times n})\, :\, {\rm MM}\, (X)=\gamma\} \end{equation}
parametrized by $\gamma\in \Gamma$.
\end{Tm}

In fact each of the sets (\ref{feb61}) contains the
element ${\rm diag}\, (1,\ldots,1,\langle \gamma, \cdot \rangle)$. In
other words, each $C(G)$-factorable matrix function can be
path-connected to an element in the above form.

Fix $\widetilde \Gamma\subset\Gamma$, a subgroup
isomorphic to $\mZ^3$ and its generators $j_1,j_2,j_3$, and define a
continuous surjective group homomorphism $\tilde p:G\to\mT^3$ by
the formula:
\[
\tilde p(g):=\bigl(\langle j_{1}, g\rangle,\langle j_{2}, g\rangle,\langle j_{3}, g\rangle\bigr)\qquad g\in G.
\]
Let $\Gamma_0$ be a finitely generated subgroup containing $\widetilde\Gamma$. Consider the family $\mathcal F$ of all finitely
generated subgroups $\Gamma_\alpha\subset \Gamma$, $\alpha\in\Lambda$, containing $\Gamma_0$, where $\Lambda$ is a suitable
index set with the partial order defined by the property that  $\alpha\leq \beta$
if and only if $\Gamma_\alpha\subseteq \Gamma_\beta$; the infimal element of $\Lambda$ is denoted
$0$ in agreement with our notation $\Gamma_0$.
We equip $\mathcal F$ with the structure of
the direct limiting system using natural inclusions
$$i_{\alpha\beta}:\Gamma_\alpha\hookrightarrow\Gamma_\beta, \qquad
\alpha\leq \beta, $$
 so that the limit
of the dual to $\mathcal F$ inverse limiting system $\widehat{\mathcal F}$ of
groups $\widehat \Gamma_\alpha$ which are Pontryagin duals
to $\Gamma_\alpha$, $\alpha\in\Lambda$, is homeomorphic to $G$. Since each
$\Gamma_\alpha$ is isomorphic to some $\mZ^{k(\alpha)}$, $k(\alpha)\in\mathbb N$, the group $\widehat\Gamma_\alpha$ is isomorphic to
the standard $k(\alpha)$-dimensional torus $\mT^{k(\alpha)}$. Thus,  $G$ is isomorphic to the inverse limit of
tori
$\mT^{k(\alpha)}$ under some continuous epimorphisms ($:=$ surjective group homomorphisms)
$$ p^{\beta}_{\alpha}:\mT^{k(\beta)}\to\mT^{k(\alpha)},
\quad k(\beta)\geq  k(\alpha), \qquad \alpha\leq \beta. $$ By $p_\alpha: G\to\mT^{k(\alpha)}$ we denote the
corresponding limit
epimorphisms.

One of the ways to describe the above construction explicitly is as follows: Let us fix a set of generators
$j_{\alpha, 1},\dots, j_{\alpha, k(\alpha)}$  in each $\Gamma_{\alpha}$. Then
\[
p_\alpha(g):=\bigl(\langle j_{\alpha, 1}, g\rangle,\dots ,\langle j_{\alpha, k(\alpha)}, g\rangle\bigr)\qquad g\in G.
\]
Now, for $\beta\ge\alpha$ and all $i$ we have: $j_{\alpha,
i}=\sum_{p=1}^{k(\beta)}c_{p,i}j_{\beta, p}$ for some
$c_{p,i}\in\mathbb Z$ (the $c_{p,i}$'s may depend on $\alpha$ and
$\beta$ as well). In particular, the homomorphism $p_\beta^\alpha$
is defined by the formula:
\[
p^\beta_\alpha(t_1,\dots, t_{k(\beta)}):=\left(\prod_{p=1}^{k(\beta)}t_p^{c_{p,1}},\dots ,
\prod_{p=1}^{k(\beta)}t_p^{c_{p, k(\alpha)}}  \right),
\]
where $(t_1,\dots, t_{k(\beta)})\in \mT^{k(\beta)}:=(\mT^1)^{k(\beta)}$
(the direct product of $k(\beta)$ circles). Also, by $\bar
p:\mT^{k(0)}\to\mT^3$ we denote similarly defined epimorphism
such that $\tilde p=\bar p\circ p_0$.

Further, the set of  pathwise connected components $C_n(G)$ of $GL(C(G)^{n\times n})$ is the direct limit of the family of
 pathwise connected
components $C_{n}(\mT^{k(\alpha)})$ of groups \penalty-10000
$GL(C(\mT^{k(\alpha)})^{n\times n})$ equipped with pullback maps
$(p^{\beta}_\alpha)^*: C_n(\mT^{k(\beta)})\to C_n(\mT^{k(\alpha)})$
(indeed, we can approximate in the uniform topology any given
element in a fixed pathwise connected component of
$GL(C(G)^{n\times n})$ by elements $X_m$, $m=1,2,\ldots$, of $P(G)$;
in turn, each $X_m$ (for sufficiently large $m$) belongs to a pathwise
connected component of some $GL(C(\mT^{\ell (m)})^{n \times n}$,
where $\ell (m)$ is the number of generators of the subgroup of
$\Gamma$  generated by the finite set $\sigma (X_m)$). These sets
have  natural group structures induced by that of $GL(\mC^{n \times
n})$. The image of the set of  pathwise connected components of
$GLF(C(G)^{n\times n})$ in the set of  pathwise connected
components of $GL(C(G)^{n\times n})$ is a subgroup generated by
pathwise connected components  containing matrix functions ${\rm
diag}(1,\dots, 1,\langle j,\cdot\rangle)$, $j\in\Gamma$ (Theorem
\ref{sep181}). In particular, this subgroup is an abelian subgroup of
$C_n(G)$ (denoted further by $C_n(G_F)$).

Assuming, on the contrary, that $GLF(C(G)^{n\times n})$ is dense in
$GL(C(G)^{n\times n})$ we obtain that $C_n(G_F)=C_n(G)$ and so
this group is abelian. Let us consider the epimorphism $\tilde p:
G\to\mT^3$. Its pullback $\tilde p^*$ induces a homomorphism
$C_n(\mT^3)\to C_n(G)$ of the corresponding groups of pathwise
connected components. Under the above assumption we obtain that
$\tilde p^*(C_n(\mT^3))\subset C_n(G_F)$. Let us take the connected
component $[X]\in C_n(\mT^3)$ containing the matrix $X:={\rm
diag}(1,\dots,1, S\circ\psi_1)$, see Example \ref{e2.4}. According to
our assumption, $\tilde p^*([X])\in C_n(G_F)$. In particular, $\tilde
p^*(X)$ belongs to the  pathwise connected component of $Y:={\rm
diag}(1,\dots, 1,\langle j,\cdot\rangle)$ for some $j\in \Gamma$.
Next, let $\Gamma_+\subset\Gamma$ be the minimal
finitely generated subgroup containing $\widetilde \Gamma$ and $j$.
Our previous direct-inverse limit construction does not change if we
start with $\Gamma_0:=\Gamma_+$.  If $\mT^{k(0)}$, $k(0)=3$ or $4$,
is the Pontryagin dual group to $\Gamma_0$, and
$p_0:G\to\mT^{k(0)}$ is the corresponding limiting map (in the
inverse limit construction), then $\tilde p^*([X])$ and $[Y]$ are
pullbacks by $p_0$ of connected components $[X']$ and $[Y']$ of $X'
=(\bar p^*(X))$, and of some $Y'={\rm diag}(1,\dots, 1,\gamma )$ for
$\gamma\in \Gamma_0$ from $C_n(\mT^{k(0)})$. Since the set of
pathwise connected components $C_n(G)$ of $GL(C(G)^{n\times n})$
is the direct limit of the family of pathwise connected components
$C_{n}(\mT^{k(\alpha)})$ of groups \penalty-10000
$GL(C(\mT^{k(\alpha)})^{n\times n})$ equipped with pullback maps
$(p^{\beta}_\alpha)^*: C_n(\mT^{k(\beta)})\to C_n(\mT^{k(\alpha)})$,
from the equality $[X]=[Y]$ it follows that there exists some
$\Gamma_\alpha$ containing $j$ such that $(\bar p\circ
p_0^\alpha)^*(X):=(p_0^\alpha)^*(X')$ belongs to the pathwise
connected component of $${\rm diag}(1,\dots,
1,(p_\alpha^0)^*(\delta)):=(p_\alpha^0)^*(Y')$$ in
$C_n(\mT^{k(\alpha)})$.

Further,
the kernel of the epimorphism $\bar p\circ p_0^\alpha:\mT^{k(\alpha)}\to\mT^3$ is a compact abelian group which is
the direct product of a torus $H_1\cong\mT^{k(\alpha)-3}$ and a
finite group $H_2$. Thus the quotient group $\mT^3(H_2):=\mT^{k(\alpha)}/H_1$ is a finite covering of $\mT^3$
with the covering group $H_2$. By $r:\mT^3(H_2) \ \rightarrow
 \mT^3$ we denote the corresponding covering map; thus, the preimage of every point in $\mT^3$ is identified with
$H_2$. We have
$$ \mT^{k(\alpha)} \ \to \ \mT^{k(\alpha)}/H_1=\mT^3(H_2) \ \stackrel{r}{\to} \ \mT^3. $$
We let also  $X_1$ be the pullback $r^*(X)$ of $X$ to $\mT^3(H_2)$.
One can easily see that $X_1$ does not belong to the  pathwise
connected component of $GL(C(\mT^3(H_2))^{n\times n})$
containing $1$. For otherwise, the map $S\circ\psi_1\circ
r:\mT^3(H_2)\to SU_2$ is homotopic to a constant map
$\mT^3(H_2)\to SU_2$. Equivalently, $\psi_1\circ
r:\mT^3(H_2)\to\mathbb S^3$ is homotopic to a constant map
$\mT^3(H_2)\to \mathbb S^3$. But the degree of $r$ is a nonzero
integer (because $r$ is a finite covering of a torus) therefore from the
latter statement and the Hopf theorem \cite[Section II.8]{Hust1} it
follows that $\psi_1:\mT^3\to\mathbb S^3$ is homotopic to a
constant map $\mT^3\to \mathbb S^3$. This contradicts our choice
of $\psi_1$, see Example \ref{e2.4}.

Next, the quotient
map $q:\mT^{k(\alpha)}\to \mT^3(H_2)$ admits a continuous
right inverse i.e. $\mT^{k(\alpha)}$ is isomorphic
(as an abelian Lie group) to $H_1\times\mT^3(H_2)$. Indeed, since all fibres of $q$ are connected, it
determines a surjective homomorphism of fundamental groups $q_*:\pi_1(\mT^{k(\alpha)})\cong \mZ^{k(\alpha)}\to
\pi_1(\mT^3(H_2))\cong\mZ^3$. In particular, there exists a monomorphism $m:\pi_1(\mT^3(H_2))\to\pi_1(\mT^{k(\alpha)})$
such that $q_*\circ m=id$ and
 $$\pi_1(\mT^{k(\alpha)})=
\pi_1(H_1)\oplus m(\pi_1(\mT^3(H_2))). $$
The map $m$ is extended by linearity to the monomorphism of $\mR$-linear hulls of the
corresponding
fundamental groups: $$\widehat m:\langle\pi_1(\mT^3(H_2))\rangle_{\mR}\cong\mR^3\to\langle\pi_1(\mT^{k(\alpha)})
\rangle_{\mR}\cong\mR^{k(\alpha)}, $$ equivariant (i.e., commuting) with respect to the actions of
$\pi_1(\mT^3(H_2))$ on $\langle\pi_1(\mT^3(H_2))\rangle_{\mR}$ and of
$\pi_1(\mT^{k(\alpha)})$ on $\langle\pi_1(\mT^{k(\alpha)})\rangle_{\mR}$
by translations. Therefore $\widehat m$ determines a homomorphism of groups
$$
\widetilde m:\mT^3(H_2)=\langle\pi_1(\mT^3(H_2))\rangle_{\mR}/\pi_1(\mT^3(H_2))\to\mT^{k(\alpha)}=
\langle\pi_1(\mT^{k(\alpha)})\rangle_{\mR}/\pi_1(\mT^{k(\alpha)}).
$$
Similarly, $\mT^{k(\alpha)}=\langle\pi_1(\mT^{k(\alpha)})\rangle_{\mR}/\pi_1(\mT^{k(\alpha)})$ and
$q$ is defined by the equivariant surjective map $\widehat q_*:
\langle\pi_1(\mT^{k(\alpha)})\rangle_{\mR}\to \langle\pi_1(\mT^3(H_2))\rangle_{\mR}$ extending $q_*$ by linearity.
From here and the definition of $q$ and $m$ we obtain that
$q\circ \widetilde m=id$.

Thus without loss of generality we may assume that $X_1$ is defined
on the subgroup $\mT^3(H_2)$ of $\mT^{k(\alpha)}$. Under this
assumption the pullback $(\bar p\circ p_0^\alpha)^*(X)$ of $X_1$ to
$\mT^{k(\alpha)}$ by means of $q$ belongs to the  pathwise
connected component containing  ${\rm diag}(1,\dots,
1,(p_\alpha^0)^*(\delta))$. Hence, $(\bar p\circ p_0^\alpha)^*(X)$
and ${\rm diag}(1,\dots, 1,(p_\alpha^0)^*(\delta))$ are homotopic
maps of $\mT^{k(\alpha)}$ into $GL(\mC^{n\times n})$. The
restriction of the homotopy between them to $\mT^3(H_2)$
determines a homotopy between ${\rm diag}(1,\dots, 1,
(p_\alpha^0)^*(\delta)|_{\mT^3(H_2)})$ and $X_1$.

But $X_1$ maps $\mT^3(H_2)$ into $SU_2$ and so its determinant is
$1$. The last two facts and Theorem \ref{sep181} imply that
$(p_\alpha^0)^*(\delta)|_{\mT^3(H_2)}$ is homotopic
to the unit homomorphism. In particular, $X_1$ belongs to the
pathwise connected component of the group
$GL(C(\mT^3(H_2))^{n\times n})$ containing $1$. This contradicts the
choice of $X$.

Thus we have proved that $p_0^*(X)$ does not belong to the set $C_n(G_F)$ and therefore
$C_n(G_F)$ is a proper subgroup of $C_n(G)$. This completes the proof of Theorem \ref{mar41}.\ \ \ \ \ $\Box$

\bigskip

We
indicate a couple of corollaries from the main results of this section:

\begin{Cy} \label{mar67}  Let $\Gamma$ and ${\cal B}$ be as in Theorem \ref{mar67'}. Then, for every natural $n\geq
2$
there
exist infinitely many pathwise connected components
of $GL({\cal B}^{n \times n})$ with the property that each one of these components does not contain
any $C(G)$-factorable element of $GL({\cal B}^{n \times n})$.
\end{Cy}

{\bf Proof.} Follows immediately from Theorem \ref{mar41} and Proposition
\ref{mar61}. $ \ \ \ \ \ \Box$
\bigskip

\begin{Cy}\label{apr181}  Let $\Gamma$ and ${\cal B}$ be as in Theorem \ref{mar67'}. Then
there exists $A\in P(G)^{n \times n}\cap GL ({\cal B}^{n \times n}) $ such that some neighborhood (in
${\cal B}^{n \times n}$) of $A$ does not contain any $C(G)$-factorable matrix function.
\end{Cy}

{\bf Proof.} Take $A_0$ in one of the pathwise connected components of
 $GL({\cal B}^{n \times n})$ described in Corollary \ref{mar67}, and approximate
$A_0$ by elements of $ P(G)^{n \times n}$. $ \ \ \ \ \ \ \Box$
\bigskip

\section{Subgroup generated by factorable \\ matrix functions}
\setcounter{equation}{0}

We focus in the section on the properties of the subgroup
generated by the factorable matrix functions. They will be needed for the proof of Theorem \ref{mar67'}, but are also of
independent interest.

An admissible algebra ${\cal B} \subseteq C(G)$ is said to be {\em
decomposing} (with respect to the order $\preceq $) if ${\cal B}={\cal
B}_+ + {\cal B}_-$. This concept was introduced in \cite{BG77}, and
studied extensively, see e.g. \cite{CG,GGK2}. Note that $W(G)$ is
decomposing while $C(G)$ is not.

\begin{La}\label{apr301} Assume ${\cal B}$ is a decomposing admissible algebra.
Then every element $X$ of the pathwise connected component
$GL_0({\cal B}^{n \times n})$ of identity in $ GL({\cal B}^{n \times n})$,
admits representation as a product of canonically ${\cal
B}$-factorable matrix functions. \end{La}

{\bf Proof.}
Let $\Lambda$
be an open neighborhood of identity in $ GL_0({\cal B}^{n \times n})$ such that every element of
 $\Lambda$ admits a canonical ${\cal B}$-factorization (existence of such $\Lambda$ follows e.g. from
\cite[Theorem XXIX.9.1]{GGK2}; this is where the decomposing
property of ${\cal B}$ is used). If $X\in  GL_0({\cal B}^{n \times n})$,
then let $X(t)$, $0\leq t\leq 1$, be a continuous path connecting $X$
with $I$ within $ GL({\cal B}^{n \times n})$; thus, $X(0)=X$, $X(1)=I$.
Partition the interval $[0,1]$, $0=t_1<t_2<\cdots <t_{p}<t_{p+1}=1$, so
that $ X_{t_j}X_{t_{j+1}}^{-1}\in \Lambda$ for $j=1,2,\ldots, p$. Then \[
X=\prod_{j=1}^p X_{t_j}X_{t_{j+1}}^{-1} \] is a product of canonically
${\cal B}$-factorable matrices. $\ \ \ \Box$
\bigskip

In the setting of decomposing admissible algebras ${\cal B}$, it is known that the
set of canonically ${\cal B}$-factorable elements of ${\cal B}^{n \times n}$ is open (this is an easy corollary of
\cite[Theorem XXIX.9.1]{GGK2} which has been used in the proof of Lemma \ref{apr301}),
and that in the classical case $G=\mT$ and ${\cal B}=W(\mT)$, the set $GLF(W(\mT)^{n \times n})$
is open, as follows from \cite{GK58}.  In general, the set $GLF(W(G)^{n \times n})$ is not open,
cf. the example at the end of Section \ref{jul51}.

Note also that for a decomposing admissible algebra ${\cal B}$, the group
${\cal G}$ that consists of all finite products of elements of $GLF({\cal B}^{n \times n})$
and their inverses, is open and closed in $GL({\cal B}^{n \times n})$. Indeed, assume $G\in GL({\cal B}^{n \times
n})$ is in the closure of ${\cal G}$; then for every $\epsilon>0$ there is $G_\epsilon \in {\cal G}$ such that
$\|G-G_\epsilon\|<\epsilon$ (the norms here are taken with respect to ${\cal B}$). Letting
$H_\epsilon:=G^{-1}G_\epsilon$ we see that $\|H_\epsilon-I\|\leq \|G^{-1}\| \epsilon$, and by taking
$\epsilon$ sufficiently small we guarantee that $H_\epsilon$ admits a canonical ${\cal B}$-factorization.
Now clearly $G=G_\epsilon H_\epsilon^{-1}$ belongs to ${\cal G}$. The openness of ${\cal G}$ can be proved by similar
arguments.
It follows that ${\cal G}$ coincides with the minimal closed subgroup of  $GL({\cal B}^{n \times n})$
that contains  $GLF({\cal B}^{n \times n})$ (again, for a decomposing admissible algebra ${\cal B}$).

For an admissible algebra ${\cal B}$, we let ${\cal D}({\cal B}^{n
\times n})$ be the union of all those pathwise connected components
of $ GL({\cal B}^{n \times n})$ that contain an element of the form
${\rm diag}\, (1,\ldots, 1, \langle j,\cdot\rangle)$, $j\in \Gamma$.
Clearly, ${\cal D}({\cal B}^{n \times n})$ is a closed and open
subgroup of $ GL({\cal B}^{n \times n})$.

\begin{Tm}\label{apr302} Let ${\cal B}$ be an admissible algebra. Then the minimal closed subgroup of $GL({\cal B}^{n
\times
n})$ containing
 $GLF ({\cal B}^{n \times n})$ is a subgroup of  ${\cal D}({\cal B}^{n \times n})$.

If in addition ${\cal B}$ is decomposing, then  the minimal closed subgroup of $GL({\cal B}^{n \times
n})$ containing
 $GLF ({\cal B}^{n \times n})$ coincides with ${\cal D}({\cal B}^{n \times n})$.
\end{Tm}

{\bf Proof.}
For the first part of the theorem it suffices to prove that every element $A\in GLF ({\cal B}^{n \times n})$
can be connected within  $GLF ({\cal B}^{n \times n})$ to an element of the form
 ${\rm diag}\, (1,\ldots, 1, \langle j,\cdot\rangle)$, $j\in \Gamma$ by a continuous path (with respect to the
${\cal B}$-norm).

Let
\begin{equation}\label{rod0''}
{A}(g)={A}_-(g)\left({\rm
diag}\,(\langle j_1,g\rangle,\ldots, \langle j_n,g \rangle)\right){A}_+(g),
\quad g\in G, \quad A_{\pm}, A_{\pm }^{-1}\in {\cal B}_{\pm}^{n \times n},
\end{equation}
be a ${\cal B}$-factorization of $A$. A known argument (see \cite[Lemma 6.4]{BrudRS} and the proof of
\cite[Theorem 6.3]{BrudRS}) shows that ${\rm
diag}\,(\langle j_1,\cdot\rangle,\ldots, \langle j_n,\cdot \rangle)$ can be connected to
${\rm
diag}\,(1,\ldots, 1, \langle j_1+\cdots +j_n,\cdot \rangle)$
within  $GLF ({\cal B}^{n \times n})$ by a continuous path.
Thus, we need only to show that
$A_\pm$ can be connected to the constant $I$ within $GL( {\cal B}_{\pm}^{n \times n})$ by a continuous path.
But this follows from Proposition \ref{mar61}, part (2), and
the proof of the first part of the theorem is completed.

For the second part,
suppose $X\in {\cal D}({\cal B}^{n \times n})$. Then
$X=X_1X_0$, where
$X_1={\rm diag}\, (1,\ldots, 1, \langle j,\cdot\rangle)$ for some $j\in \Gamma$ and
$X_0\in  GL_0({\cal B}^{n \times n})$. Now $X_1$ is obviously ${\cal B}$-factorable, and
$X_0$ is a product of canonically ${\cal B}$-factorable matrices by Lemma \ref{apr301}. Thus,
${\cal D}({\cal B}^{n \times n})$ is contained in the  minimal closed subgroup of $GL({\cal B}^{n \times n})$ containing
all ${\cal B}$-factorable matrices.
$\ \ \ \Box$
\bigskip

Notice that Lemma \ref{apr301} and Theorem \ref{apr302} are valid for any
connected compact abelian group $G$, not only for those whose dual
group $\Gamma$ contains ${\mathbb Z}^3$.

The case of the two-dimensional torus is perhaps interesting:
\begin{Cy}
Let $G=\mT^2$, and let ${\cal B}$ be a decomposing admissible algebra.
Then  $GL({\cal B}^{n \times
n})$ is equal to
the minimal closed subgroup of $GL({\cal B}^{n \times
n})$ containing
 $GLF ({\cal B}^{n \times n})$.
\end{Cy}

{\bf Proof.}
By Theorem \ref{apr302}, the minimal closed subgroup of $GL({\cal B}^{n \times
n})$ containing
 $GLF ({\cal B}^{n \times n})$ coincides with ${\cal D}({\cal B}^{n \times n})$. Using
Proposition \ref{mar61} and the proof of Proposition \ref{dec183} we see that
in turn ${\cal D}({\cal B}^{n \times n})$ coincides with
 $GL({\cal B}^{n \times
n})$. Indeed, the pathwise connected components of both groups are
parametrized identically by
${\mathbb Z}^2$.
$\ \ \ \Box$
\bigskip

{\bf Proof of Theorem \ref{mar67'}}. As we have seen in the proof of Theorem \ref{mar41} (see also
Proposition \ref{mar61}), there
are infinitely many
pathwise connected components of $GL({\cal B}^{n\times n})$ that do not intersect  ${\cal D}({\cal B}^{n \times n})$.
The result now follows from Theorem \ref{apr302}.
$\ \ \ \Box$

\section{Scalar valued functions}
\setcounter{equation}{0}

In this section we consider the scalar case.

\begin{Tm}\label{may111} Let ${\cal B}\subseteq C(G)$ be an admissible algebra, where $G$ is a connected compact abelian group.
Then:
\begin{itemize}
\item[${\rm (a)}$]
The set
$GLF({\cal B})$ of ${\cal B}$-factorable scalar functions is dense in $GL({\cal B})$;
\item[${\rm (b)}$]
The equality $GLF({\cal B})=GL({\cal B})$ holds if and only if ${\cal B}$ is decomposing.
\end{itemize}
\end{Tm}

In the classical case $G={\mathbb T}$ part (b) is well known, see e.g. \cite[Theorem 3.1]{GKru92}.

For the proof of Theorem \ref{may111} it will be convenient to start with preliminary results. The next proposition is
perhaps independently interesting and holds for matrix functions as well.

\begin{Pn}\label{may172}  Let ${\cal B}$ be an admissible algebra,  $f\in {\cal B}^{n \times n}$, and let
$$
\Omega = \{z \in {\mathbb C}\ :\  z \mbox{ is an eigenvalue of $f(g)$ for
some }
g \in G\}.
$$
If $\Psi$ is an analytic function in an open neighborhood of the
closure of $\Omega$, then $\Psi\circ f \in {\cal B}^{n \times n}$.
\end{Pn}

Here, for every fixed $g\in G$, $\Psi\circ f(g)$ is understood as
the
$n \times n$ matrix defined by the standard functional calculus.

The proof is essentially the same as that of \cite[Proposition 2.3]{RSW01ot}, \cite[Proposition 2.3]{RSW98}, where
it was proved for $W(G)$, and with $\Gamma$ an additive subgroup of ${\mathbb R}^k$ (with the discrete topology).

\begin{Pn}\label{may191} A character $c\in \Gamma$ is an exponential, i.e.  $\langle c,\cdot\rangle =e^{f(\cdot)}$ for some
$c\in \Gamma$ and $f\in C(G)$, if and only if $c=0$.
\end{Pn}

{\bf Proof.} For the reader's convenience we supply a known proof. The ``if" part is obvious.
For the ``only if" part, we have
$$
e^{f(gh)}=<c,gh>=<c,g>\cdot <c,h>=e^{f(g)+f(h)}, \qquad \forall \ \ g,h\in G. $$
Therefore
\begin{equation}\label{may194}
f(gh)-f(g)-f(h)=2\pi i a(g,h),
\end{equation}
where $a(g,h)$ is an integer valued  continuous function on $G\times G$. Since $G$ is connected, $G\times G$ is also connected
and
therefore $a$ is a constant.
Now if $g=h=1$, then (\ref{may194}) gives
$$
a(g,h)=a(1,1)=\frac{-f(1)}{2\pi i}.
$$
For a fixed $g\in G$ by applying inductively (\ref{may194}) we have for any natural $n$:
\begin{equation}\label{may195}
f(g^n)= nf(g)-(n-1)f(1)=f(g)+(n-1)(f(g)-f(1)).
\end{equation}
Since $G$ is compact, the image $f(G)$ is a bounded subset of the plane. On
the other hand (\ref{may195}) shows that this is possible for $n\ \rightarrow \ \infty$ if and
only if $f(g)=f(1)$ for any $g\in G$. Thus $f$ is constant and moreover $f\in 2\pi i {\mathbb Z}$.
 From here we get $<c,\cdot>=e^{f(1)}=1.$ $\ \ \ \Box$
\bigskip

{\bf Proof of Theorem \ref{may111}.} Part (a). We use approximation
(in the norm of ${\cal B}$) of any given element of $GL({\cal B})$ by
elements of $P(G)\cap GL({\cal B})$, thereby reducing the proof to the
case of finitely generated subgroups $\Gamma'$ of $\Gamma$ rather
than $\Gamma$ itself.

Indeed, letting  $C(G)_{\Gamma'}$ stand for the closure in $C(G)$ of
the set $P(G)_{\Gamma'}$ of elements of  $P(G)$ with the Bohr-Fourier
spectra in $\Gamma'$, we first note that  $C(G)_{\Gamma'}$ coincides
with the set of elements of $C(G)$ having the Bohr-Fourier spectra in
$\Gamma'$ (see e.g. \cite[Theorem 7.14]{C}). Next, let  $G'$ be the
dual group of $\Gamma'$, and  observe the identification  $C(G')\cong
C(G)_{\Gamma'}$. This follows easily from the fact that $\Gamma'$ is
the dual group of $G'\cong G/H$, where $H:=\{g\in G\colon\langle
\gamma, g\rangle =1 \quad \forall \ \ \gamma \in \Gamma'\}$ is the
annihilator of $\Gamma'$ and a closed subgroup of $G$ \cite[Section
2.1]{Ru90}; thus, if $a\in P(G)_{\Gamma'}$, then $a$ is constant on
every coset of $G$ by $H$, and therefore every element of
$C(G)_{\Gamma'}$ is also constant on every coset of $G$ by $H$.
Also, using the approximation of any given element of $GL({\cal
B})$ by elements of $P(G)\cap GL({\cal B})$, we replace ${\cal B}$
with the closure (in the norm of ${\cal B}$) of $P(G)_{\Gamma'}$;
denote this closure by ${\cal B}'$. Note that ${\cal B}'$ can be understood also as a
subalgebra of $C(\Gamma')$ via the above identification
 $C(G')\cong C(G)_{\Gamma'}$; the inverse closedness of ${\cal B}'$ follows by identifying
the elements of ${\cal B}'$ with those  functions in ${\cal B}$ that are constant on every coset of $G$ by $G'$.
The
finitely generated subgroups $\Gamma'$ are further
identified with $\mZ^q$, for suitable integers $q$.

Next,
we prove that every element $a\in P(\mT^q) \cap GL({\cal B})$ admits a
decomposition of the form
\begin{equation}\label{may171}
a=|a|\langle c, \cdot \rangle e^{u}=e^b\langle c, \cdot \rangle e^{u},
 \end{equation}
for some $c\in \mZ^q$ and
$u,b\in {\cal B}$. The proof of (\ref{may171}) will mimic that of
 \cite[Proposition 2.4]{RSW03}, \cite[Theorem 2.4]{RSW01ot}, where (\ref{may171}) was proved
for the case ${\cal B}= W(\mT^q)$.
First note that $a\overline{a}\in  {\cal B}$, hence by Proposition \ref{may172},
$|a|=\sqrt{a\overline{a}}$ and $\log |a|$ belong to $ {\cal B}$. Thus, the second equality in
(\ref{may171})
follows. For the first equality, note that it was proved in \cite{perov} with $u\in C(\mT^q)$.
Thus, we need only to show that in fact $u\in  {\cal B}$.
Represent the function $u$ as the sum $u=u_0+u_1$, where
$u_0\in  {\cal B}$ and $\|u_1\|_{C(\mT^q)}<\frac{\pi}{2}$.
Then
\begin{equation}\label{bes1}
e^{u_1(g)}=a(g)|a(g)|^{-1}\langle -c, g\rangle e^{-u_0(g)}, \qquad g\in \mT^q.
\end{equation}
By Proposition \ref{may172}
the function $|a|^{-1}$ belongs to $ {\cal B}$.
The other two multiples in the right hand side of (\ref{bes1}) obviously
belong to $ {\cal B}$; hence, so does the function $\xi=e^{u_1}$.
On the other hand, the values of $\xi$ lie in the right open half-plane.
Using Proposition \ref{may172} again, we may define $u_2\in  {\cal B}$ so that
$e^{u_2}=\xi$. Since $u_1$ and $u_2$ are both continuous on $\mT^q$, this
means that they differ by a constant summand. But then $u_1$ belongs to
$ {\cal B}$ simultaneously with $u_2$. Finally, the function $u=u_0+u_1$
belongs to $ {\cal B}$ as well.

Now, for  $a\in P(\mT^q) \cap GL({\cal B})$ and its decomposition (\ref{may171}),
let $\{u_m\}_{m=1}^\infty$, $\{b_m\}_{m=1}^\infty$ be sequences such that
$u_m,b_m\in {\cal B}_++{\cal B}_-$ (where ${\cal B}_++{\cal B}_-$ is considered with its Banach
algebra norm) and
$$ \lim_{m\rightarrow \infty} u_m =u, \qquad \lim_{m\rightarrow \infty} b_m =b $$
in the ${\cal B}$-norm. Then clearly
\begin{equation}\label{may171'} a_m:=e^{b_m}\langle c, \cdot \rangle e^{u_m} \quad \longrightarrow \quad a
\end{equation} in the  ${\cal B}$-norm
as $m\rightarrow\infty$. On the other hand, it follows from
(\ref{may171'}) that every  $a_m$ admits a ${\cal B}$-factorization.
Indeed, write $u_m=u_{m+}+u_{m-}$, $b_m=b_{m+}+b_{m-}$ with
$u_{m\pm }, b_{m\pm}\in {\cal B}_\pm$. Then
$$ a_m=e^{b_{m-}}e^{u_{m-}}\langle c, \cdot \rangle e^{b_{m+}}e^{u_{m+}} $$
is the desired factorization. This proves the density of $GLF({\cal B})$ in $GL({\cal B})$.
\bigskip

Part (b), the ``if" part. Assume ${\cal B}$ is decomposing, and let $a\in
GL({\cal B})$. Approximate $a$ in the ${\cal B}$-norm by a function
$b\in P(G)\cap  GL({\cal B})$ for which the function $ab^{-1}$ is so
close to identity that it is canonically ${\cal B}$-factorable (see e.g.
\cite[Theorem XXIX.9.1]{GGK2}; the decomposing property of ${\cal
B}$ is essential here). As proved in the preceding paragraph, we can
approximate $b$ by $c\in GL({\cal B})$, again in the ${\cal B}$-norm,
so that $c$ admits a ${\cal B}$-factorization and $bc^{-1}$ is
canonically
 ${\cal B}$-factorable (we use \cite[Theorem XXIX.9.1]{GGK2} again).
Now $a=c\cdot (bc^{-1})\cdot (ab^{-1})$, and since every factor admits a ${\cal B}$-factorization,
the ${\cal B}$-factorability of $a$ follows.
\bigskip

Part (b), the ``only if" part. Assume that every invertible element of
${\cal B}$ is ${\cal B}$-factorable.  Let
$h\in {\cal B}$ and let $f=e^h$. Since $f\in GL({\cal B})$, we have a
${\cal B}$-factorization $ f=b_{-}\langle c,\cdot\rangle  b_{+}$, $c\in
\Gamma$, $b_\pm \in GL({\cal B}_\pm)$. By  Proposition~\ref{mar61}, part (2),
$b_\pm$ belongs to the pathwise connected component of identity in
$GL({\cal B}_\pm)$. Thus, the character $\langle c,\cdot\rangle $
belongs to the same pathwise connected component of $GL({\cal
B})$ as $e^h$ does. But this pathwise connected component is just
the component that contains the constant $1$, or equivalently, the set
of all exponentials $e^s$, $s\in {\cal B}$. Thus, $\langle c,\cdot\rangle
$ is an exponential, and by Proposition~\ref{may191} we obtain $c=0$.
Furthermore, since $b_\pm$ belongs to the pathwise connected
component of identity in $GL({\cal B}_\pm)$, we have
$b_\pm=e^{d_\pm}$ for some $d\in {\cal B}_\pm$. Now
$e^h=e^{d_-}e^{d_+}$, and so $h=d_-+d_++2\pi i k$, where $k\,:\, G \,
\rightarrow {\mathbb Z}$ is a continuous function. So $k$ must be
constant, and we are done. $\ \ \ \Box$

\section{Non-denseness of triangularizable matrix functions}
\setcounter{equation}{0}

In this section we apply our main result to the problem of
triangularizability. As everywhere in the paper, let $G$ be a compact
abelian group with the dual group $\Gamma$. Let ${\cal B}$ be an
admissible algebra. An element $A\in {\cal B}^{n \times n}$ is said to
be (right) ${\cal B}$-{\em triangularizable} if $A$ admits a
representation (\ref{rod0}), where the diagonal term ${\rm
diag}\,(\langle j_1,g\rangle,\ldots, \langle j_n,g\rangle )$ is replaced
by a triangular matrix $T=[t_{ij}]_{i,j=1}^n$, with $t_{ij}\in {\cal B}$ for
$i,j=1,\ldots n$, $t_{ij}=0$ if $i>j$, and the diagonal elements $t_{11},
\ldots, t_{nn}$ belong to $GL({\cal B})$. Denote by $GLT({\cal B}^{n
\times n})$ the set of $n\times n$  ${\cal B}$-triangularizable matrix
functions. Clearly,
$$ GLF({\cal B}^{n \times n})\subseteq GLT({\cal B}^{n \times n})\subseteq GL({\cal B}^{n \times n}). $$
An example (\cite{KarlSpit85}, see also \cite[Theorem 8.17]{BKS1}) shows that, for
$\Gamma={\mathbb R}$, there exist $2 \times 2$ triangular matrix functions in $P(G)\cap GL(W(G)^{2 \times 2})$
that are not $W(G)$-factorable.

The following question has been proposed, in the context of the algebra $W(G)$: {\em Does $ GL({\cal B}^{n \times
n})= GLT({\cal B}^{n \times n})$ hold for admissible algebras? } Theorem \ref{mar67'} implies that generally the answer is no.
More precisely, we have the following result.
Denote by  ${\cal T}({\cal B}^{n \times n})$ the minimal closed subgroup of $GL({\cal B}^{n \times n})$ that contains
 $GLT ({\cal B}^{n \times n})$.

\begin{Tm} \label{mar67''}  Let $\Gamma$ be a torsion free abelian group that contains a subgroup
isomorphic to
${\mathbb Z}^3$, and let ${\cal B}$ be an admissible algebra. Then,
for every natural $n\geq 2$
there
exist infinitely many pathwise connected components
of $GL({\cal B}^{n \times n})$ with the property that each one of these components does not
intersect
 ${\cal T}({\cal B}^{n \times n})$.
In particular,  ${\cal T}({\cal B}^{n \times n})$ is not dense in  $GL({\cal B}^{n \times n})$.
 \end{Tm}

{\bf Proof.}
Let $A\in GLT({\cal B}^{n \times n})$, so
$$ A(g)={A}_-(g)\, \left([t_{ij}(g)]_{i,j=1}^n\right)\, {A}_+(g),
\quad g\in G, $$
where $A_{\pm}, A_{\pm }^{-1}\in {\cal B}_{\pm}^{n \times n}$,
 $t_{ij}\in {\cal B}$, $t_{ij}=0$ if $i>j$, and $t_{jj}\in GL({\cal B})$.
By the pathwise connectivity of  $GL({\cal B}_{\pm}^{n \times n})$
(Proposition~\ref{mar61} part (2)), we see that there is a continuous
path in $GL({\cal B}^{n \times n})$ connecting $A$ with
$[t_{ij}]_{i,j=1}^n$. In turn, $[t_{ij}]_{i,j=1}^n$ is obviously continuously
connected in $GL({\cal B}^{n \times n})$ with ${\rm diag}\,
(t_{11},\ldots, t_{nn})$. By Theorem \ref{may111}, each $t_{jj}$ can be
approximated arbitrarily well (in the ${\cal B}$-norm) by ${\cal
B}$-factorable scalar functions, hence ${\rm diag}\, (t_{11},\ldots,
t_{nn})$ can be approximated arbitrarily well by ${\cal B}$-factorable
matrix functions. Thus, $A$ belongs to a pathwise connected
component of $GL({\cal B}^{n \times n})$ that contains a ${\cal
B}$-factorable matrix function. By Theorem \ref{apr302}, the union of
those pathwise connected components of $GL({\cal B}^{n \times n})$
is contained in ${\cal D}({\cal B}^{n \times n})$, and since  ${\cal
D}({\cal B}^{n \times n})$ is an open and closed subgroup of
$GL({\cal B}^{n \times n})$, also  ${\cal T}({\cal B}^{n \times n})$ is
contained in
 ${\cal D}({\cal B}^{n \times n})$. It remains to observe that
there are infinitely many pathwise connected components of
$GL({\cal B}^{n\times n})$ that do not intersect  ${\cal D}({\cal B}^{n
\times n})$ (cf. the proof of Theorem \ref{mar67'}). $\ \ \ \Box$

\bibliographystyle{amsplain}
\providecommand{\bysame}{\leavevmode\hbox to3em{\hrulefill}\thinspace}
\providecommand{\MR}{\relax\ifhmode\unskip\space\fi MR }
\providecommand{\MRhref}[2]{%
  \href{http://www.ams.org/mathscinet-getitem?mr=#1}{#2}
}
\providecommand{\href}[2]{#2}

\end{document}